\documentclass[11pt]{article}

\usepackage{amsfonts,amsmath,amsthm,latexsym,color,epsfig,mathrsfs}
\newtheorem{thm}{Theorem}

\title{On the normalized Shannon capacity of a union}
\author{Peter Keevash
\thanks{Mathematical Institute, University of Oxford, Oxford, UK. Email: keevash@maths.ox.ac.uk.
Research supported in part by ERC Consolidator Grant 647678.}
\and Eoin Long \thanks{School of Mathematical Sciences, Tel Aviv University, Tel Aviv, Israel. Email: eoinlong@post.tau.ac.il.}}

\begin{document}
\maketitle

\begin{abstract}
Let $G_1 \times G_2$ denote the strong product of graphs $G_1$ and $G_2$,
i.e.\ the graph on $V(G_1) \times V(G_2)$ in which $(u_1,u_2)$ and $(v_1,v_2)$
are adjacent if for each $i=1,2$ we have $u_i=v_i$ or $u_iv_i \in E(G_i)$.
The Shannon capacity of $G$ is $c(G) = \lim_{n\to \infty} \alpha (G^n)^{1/n}$,
where $G^n$ denotes the $n$-fold strong power of $G$,
and $\alpha (H)$ denotes the independence number of a graph $H$.
The normalized Shannon capacity of $G$ is $C(G) = \frac {\log c(G)}{\log |V(G)|}$. 
Alon \cite{alon} asked whether for every $\epsilon > 0$ there are graphs $G$ and $G'$ 
satisfying $C(G), C(G') < \epsilon$ but with $C(G + G') > 1 - \epsilon $.
We show that the answer is no.
\end{abstract}

Despite much impressive work (e.g.\ \cite{alon}, \cite{alon lubetzky}, \cite{alon orlitsky},  \cite{haemers}, \cite{lovasz})
since the introduction of the Shannon capacity in \cite{shannon}, many natural questions regarding this parameter
remain widely open (see \cite{graph powers}, \cite{korner orlitsky} for surveys).
Let $G_1 + G_2$ denote the disjoint union of the graphs $G_1$ and $G_2$. 
It is easy to see that $c(G_1 + G_2) \geq c(G_1) + c(G_2)$. 
Shannon \cite{shannon} conjectured that $c(G_1 + G_2) = c(G_1) + c(G_2)$,
but this was disproved in a strong form by Alon \cite{alon} who showed
that there are $n$-vertex graphs $G_1,G_2$ 
with $c(G_i) < e^{c\sqrt{\log n \log\log n}}$ but $c(G_1 + G_2) \geq \sqrt n$.
In terms of the normalized Shannon capacity, 
this implies that for any $\epsilon >0$, there exist graphs $G_1$, $G_2$ 
with $C(G_i) < \epsilon $ but $C(G_1 + G_2) > 1/2 -\epsilon$. 
Alon \cite{alon} asked whether `$1/2$' can be changed to `$1$' here.
In this short note we will give a negative answer to this question.
In fact, the following result implies that `$1/2$' is tight.

\begin{thm}
\label{normalized shannon capacity}
If $C(G_1) \le \epsilon$ and $C(G_2) \le \epsilon$
then $C(G_1 + G_2) \leq \frac{1 + \epsilon }{2} + \frac{1-\epsilon}{2\log _2 (|V(G_1)| + |V(G_2)|)}$.
\end{thm}

\noindent \textbf{Proof.}
Let $N_i = |V(G_i)|$ for $i=1,2$.
Fix a maximum size independent set $I$ in $(G_1+G_2)^n$ for some $n\in {\mathbb N}$.
We write $|I|=\sum_{S \subset [n]} |I_S|$, where
$I_S = \{ x=(x_1,\ldots,x_n) \in I: x_i \in V(G_1) \Leftrightarrow i \in S\}$.

To bound $|I_S|$, we may suppose that $S=[m]$ for some $0\leq m \leq n$.
Then $I_S$ is an independent set in $G_1^{m} \times G_2^{n-m}$. 
As $C(G_1) \leq \epsilon $, by supermultiplicativity $\alpha (G_1^{m}) \leq N_1^{\epsilon m}$;
similarly, $\alpha (G_2^{n-m}) \leq N_2^{\epsilon (n-m)}$.
For any $x \in V(G_1)^m$, the set of $y \in V(G_2)^{n-m}$ 
such that $(x,y) \in I_S$ is independent in $G_2^{n-m}$,
so $|I_S| \le N_1^m N_2^{\epsilon (n-m)}$.
Similarly, $|I_S| \le N_1^{\epsilon m}N_2^{n-m}$.

We multiply these bounds:
$|I_S|^2 \le (N_1^m N_2^{n-m})^{1+\epsilon}$.
Writing $\gamma = \frac{N_1}{N_1 + N_2}$, we have
\begin{eqnarray*}
\alpha ((G_1 + G_2)^n) = |I| & = & \sum _{S \subset [n]} |I_S| 
\leq \sum _{m=0}^n \binom {n}{m} \left(N_1^{(1 +\epsilon )/2}\right)^m \left(N_2^{(1 +\epsilon )/2}\right)^{n-m}\\ 
& = & (N_1^{(1 +\epsilon )/2} + N_2^{(1 +\epsilon )/2})^{n} \\
& = & (\gamma ^{(1 +\epsilon )/2} + (1-\gamma )^{(1 +\epsilon )/2})^n(N_1 + N_2)^{(1+\epsilon)n/2}\\
& \leq & 2^{(1-\epsilon )n/2} (N_1 + N_2)^{(1+\epsilon)n/2},
\end{eqnarray*}
as $\gamma ^{b} + (1-\gamma )^{b}$ is maximized at $\gamma = 1/2$ for $0<b<1$ and $0 \leq \gamma \leq 1$. 
Therefore
\begin{equation*}
C(G_1 + G_2) = \lim _{n\to \infty} \frac{\log \alpha ((G_1 + G_2 )^{n})}{n \log {(N_1 + N_2)}} \leq \frac{1 + \epsilon }{2} + \frac{1 - \epsilon }{2\log _2(N_1 + N_2)}.
\end{equation*}

\end{document}